\documentclass[10pt]{amsart}
\usepackage{latexsym, amsmath,amssymb}
\usepackage{graphicx}
\usepackage{enumerate}
\usepackage{color}

\renewcommand{\epsilon}{\varepsilon}

\newtheorem{theorem}{Theorem}[section]

\newtheorem{lemma}[theorem]{Lemma}

\newtheorem{corr}[theorem]{Corollary}
\newtheorem{prop}[theorem]{Proposition}
\newtheorem{deff}[theorem]{Definition}

\newcommand{\bth}{\begin{theorem}}
	\newcommand{\ble}{\begin{lemma}}
		\newcommand{\bcor}{\begin{corr}}

			\newcommand{\bdeff}{\begin{deff}}

				\newcommand{\bprop}{\begin{proposition}}
					\newcommand{\ele}{\end{lemma}}
				\newcommand{\ecor}{\end{corr}}
			\newcommand{\edeff}{\end{deff}}
		
		\newcommand{\eprop}{\end{proposition}}

	\newcommand{\la}{\lambda}
	
	\newcommand{\eps}{\varepsilon}

	\newcommand{\supp}{\text{supp }}
	\renewcommand{\Pi}{\varPi}

	\renewcommand{\epsilon}{\varepsilon}

	\newcommand{\ls}{\lesssim}
	
	\newcommand{\1}{{\rm 1\hspace*{-0.4ex}%
			\rule{0.1ex}{1.52ex}\hspace*{0.2ex}}}
	
	\newcommand{\ola}{\1_\la}
	
	\newcommand{\ala}{\tilde \1_\la}

	\numberwithin{equation}{section}

\begin{document}
		\title[]{Weyl laws for Schr\"odinger operators on compact manifolds with boundary}
		\author{Xiaoqi Huang, Xing Wang and Cheng Zhang}
		\address{Department of Mathematics\\
			Louisiana State University\\
			Baton Rouge, LA 70803, USA}
		\email{xhuang49@lsu.edu}
		\address{Department of Mathematics\\
			Hunan University\\
			Changsha, HN 410012, China}
		\email{xingwang@hnu.edu.cn}
		\address{Mathematical Sciences Center\\
			Tsinghua University\\
			Beijing, BJ 100084, China}
		\email{czhang98@tsinghua.edu.cn}
		\date{}
		\begin{abstract} We prove Weyl laws for Schr\"odinger operators with critically singular potentials  on compact manifolds with boundary. We also improve the Weyl remainder estimates under the condition that the set of all periodic geodesic billiards has measure 0. These extend the classical results by Seeley \cite{seeley1, seeley2}, Ivrii \cite{ivrii} and Melrose \cite{mel}. The proof uses the Gaussian heat kernel bounds for short times and a perturbation argument involving the wave equation.
		\end{abstract}
		
		\maketitle

\section{Introduction}
	Let $(M,g)$ be a smooth compact Riemannian manifold of dimension $n\ge2$ with smooth boundary $\partial M$. Let $\Delta_g$ be the Laplace-Beltrami operator on $M$.  Let $\nu$ be the outward unit normal vector field along $\partial\Omega$. Under either Dirichlet ($u|_{\partial\Omega}=0$), Neumann ($\partial_\nu u|_{\partial\Omega}=0$) or Robin ($(\partial_\nu u+\sigma u)|_{\partial\Omega}=0$ with nonnegative $\sigma\in C^\infty(\partial\Omega)$) boundary condition, the Laplacian $-\Delta_g$ is self-adjoint and nonnegative on its domain, and has discrete spectrum $\{\la_j\}_{j=1}^\infty$, 	where the eigenvalues, $\la_1\le \la_2\le \cdots$, are arranged in increasing order and we account for multiplicity.  See e.g. Taylor \cite{taylorpde1}. 
	
	Under either Dirichlet, Neumann or Robin boundary condition, the Weyl law  for the Laplacian $-\Delta_g$   is the following one-term asymptotic formula
	\begin{equation}\label{weyllaw}\# \{j: \, \la_j\le \la\}=(2\pi)^{-n}\omega_n |M|\lambda^n+O(\lambda^{n-1}),\end{equation}
	where $\omega_n$ is the volume of the unit ball in $\mathbb{R}^n$ and $|M|$ is the  Riemannian volume of $M$. The formula \eqref{weyllaw}  with the sharp remainder term $O(\la^{n-1})$ is due to Seeley \cite{seeley1,seeley2}. Indeed, he constructed a short-time parametrix for the wave equation near the boundary under either Dirichlet or Neumann boundary condition, and proved \eqref{weyllaw} by a Tauberian argument. Under Robin boundary condition, the eigenvalues $\la_j^R$ lie between the Dirichlet eigenvalues $\la_j^D$ and the Neumann eigenvalues $\la_j^N$. So the formula \eqref{weyllaw} remains valid. See also Weyl \cite{weyl}, Courant \cite{courant}, Carleman \cite{car1,car2}, Avakumovi\'c \cite{ava}, Levitan \cite{lev}, H\"ormander \cite{horm}, B\'erard \cite{berard} and many others for related works on compact manifolds with or without boundary. An extensive bibliographical review can be found in \cite{RSS}.
	
	Weyl \cite{weylconj} put forward a conjecture on the following two-term asymptotic formula	\begin{equation}\label{weylconj}	\#\{j:\la_j\le \la\}=(2\pi)^{-n}\omega_n|M|\la^n\mp\frac14(2\pi)^{1-n}\omega_{n-1}|\partial M|\la^{n-1}+o(\la^{n-1}),\end{equation}
	where the minus corresponds to the Dirichlet condition, the plus to the Neumann condition, and $|\partial M|$ is the $(n-1)$-dimensional volume of $\partial M$. Note that \eqref{weylconj} cannot hold on the standard sphere or hemisphere since the eigenvalues have very high multiplicities. However, it is still open for general bounded domains in $\mathbb{R}^n$.  Duistermaat-Guillemin \cite{DG} proved \eqref{weylconj} on closed manifolds $(\partial M=\emptyset)$, under the assumption that 
	\begin{equation}\label{DGcond}\text{the set of all periodic geodesics has measure 0.}\end{equation} 
	 Ivrii \cite{ivrii} proved \eqref{weylconj} on compact manifolds with boundary, under the assumption that 
	\begin{equation}\label{ivrcond}\text{the set of all periodic geodesic billiards has measure 0,}\end{equation} which generalizes the condition \eqref{DGcond}.  The condition \eqref{ivrcond} is only known to be true for special classes of shapes of domains in $\mathbb{R}^n$, such as the convex domains with analytic boundary. If one allows piecewise smooth boundaries, then \eqref{ivrcond} is also true if each smooth component of $\partial M$ has nonpositive normal curvature (say, a polyhedron satisfies this condition). See Safarov-Vassiliev \cite[Theorem 1.6.1 \& Conjecture 1.3.35]{safa} and references therein. See also Melrose \cite{mel}, H\"ormander \cite[Corollary 29.3.4]{hor4}, Ivrii \cite{ivrii2}.
	
	In this paper, we shall extend the classical results \eqref{weyllaw} and \eqref{weylconj} to Schr\"odinger operators with critically singular potentials.

	\subsection{Schr\"odinger operators}We shall assume throughout that the potentials $V$ are real-valued and 
	\begin{equation}\label{condV}
		V\in L^1(M)\ \ \text{with}\ \  V^-\in \mathcal{K}(M).
	\end{equation} Here $V^-=\max\{0,-V\}$ and  $\mathcal{K}(M)$ is the Kato class. Recall that $\mathcal{K}(M)$ is all $V$ satisfying  
	\begin{equation}\label{kato}\lim_{\delta\to 0}\sup_{x\in M}\int_{d_g(y,x)<\delta}|V(y)|W_n(d_g(x,y))dy=0,\end{equation}
	where \[W_n(r)=\begin{cases}r^{2-n},\quad\quad\quad\quad n\ge3\\
		\log(2+r^{-1}),\ \ n=2\end{cases}\]
	and $d_g$, $dy$ denote the geodesic distance and the volume element on $(M,g)$, respectively.
	By H\"older inequality, we have $L^{q}(M)\subset \mathcal{K}(M)\subset L^1(M)$ for all $q>\frac n2$. So the condition \eqref{condV} is weaker than $V\in \mathcal{K}(M)$.
	
	It is known that the condition \eqref{condV} can ensure the Schr\"odinger operator $$H_V=-\Delta_g+V$$ is self-adjoint and bounded from below. Since $M$ is compact, the spectrum of $H_V$ is discrete, and the associated eigenfunctions are bounded. If $V\in\mathcal{K}(M)$, then the eigenfunctions are also continuous. These results rely on the Gaussian heat kernel bound \eqref{gauss} for short times. See \cite{hs}, \cite{BSS}, \cite{simonsurvey}, \cite{LHB}, \cite{Gu17}. 
	
	 After adding a constant to the potential we may, and always shall assume that $H_V$ is bounded from below by one. This just shifts the spectrum and does not change the eigenfunctions. We shall write the spectrum
	of $\sqrt{H_V}$  as
		$$\{\tau_k\}_{k=1}^\infty,$$
	where the eigenvalues, $\tau_1\le \tau_2\le \cdots$, are arranged in increasing order and we account for multiplicity.  For each $\tau_k$ there is an 
	eigenfunction $e_{\tau_k}\in \text{Dom }(H_V)$ (the domain of $H_V$) so that
	\begin{equation}\label{1.6}
		H_Ve_{\tau_k}=\tau^2_ke_{\tau_k},\ \ {\rm and}\  \ \int_M |e_{\tau_k}(x)|^2 \, dx=1.
	\end{equation} 
	To  be consistent, we shall let
$$	H^0=-\Delta_g$$
be the unperturbed operator.  The corresponding eigenvalues and associated $L^2$-normalized
eigenfunctions are denoted by $\{\lambda_j\}_{j=1}^\infty$ and $\{e^0_j\}_{j=1}^\infty$, respectively so
that
\begin{equation}\label{1.9}
	H^0e^0_j=\lambda^2_j e^0_j, \quad \text{and }\, \, 
	\int_M |e^0_j(x)|^2 \, dx=1.
\end{equation}
Both $\{e_{\tau_k}\}_{k=1}^\infty$ and $\{e^0_j\}_{j=1}^\infty$ are orthonormal bases for $L^2(M)$.  Let $P^0=\sqrt{H^0}$ and $P_V=\sqrt{H_V}$.  

\subsection{Heat kernel bounds}
In this paper, we shall only use the heat kernel bounds for short times. Under either Dirichlet, Neumann or Robin boundary condition,  we have
\begin{equation}\label{gauss}|e^{t\Delta_g}(x,y)|\ls t^{-\frac n2}e^{-cd_g(x,y)^2/t},\ \ 0<t\le 1,\end{equation}
for some constant $c>0$. These  Gaussian heat kernel bounds for short times were proved by Greiner \cite{Gre}. See also Li-Yau \cite{LY}, Davies \cite{davies}, Daners \cite{dan}.  In \cite{Gre}, Greiner constructed the parametrix of the heat equation, and used it to calculate the asymptotic expansion for the heat trace as $t\to 0^+$. This is the approach exploited by  McKean-Singer \cite{MS} to solve  Kac's conjecture \cite{kac} on  the heat trace.  Moreover, on physical grounds one expects that for short times the heat kernel is dominated
by local contributions that do not involve the boundary. This is essentially the
principle of not feeling the boundary by Kac \cite{kac2}.  

By mimicking the proof of the Feynman-Kac formula \cite{Simonbook} and 
the Gaussian heat kernel bound \cite[Prop. B.6.7] {simonsurvey} in $\mathbb{R}^n$, one may obtain the following short-time bound for the Schr\"odinger heat kernel under the condition \eqref{condV} on manifolds with or without boundary. See e.g. Sturm \cite{sturm}.
\begin{lemma}\label{heatk}
	Let $V\in L^1(M)$ with $V^-\in \mathcal{K}(M)$. Under either Dirichlet, Neumann or Robin boundary condition, the  Schr\"odinger heat kernel satisfies the Gaussian upper bound
	\begin{equation}\label{heatV}
		|e^{-tH_V}(x,y)|\ls t^{-\frac n2}e^{-c_1d_g(x,y)^2/t},\ \ 0<t\le 1,
	\end{equation}
	for some constant $c_1>0$.
\end{lemma}
By \eqref{heatV}, we have the eigenfunction bound \begin{equation}\label{rough}
	\sum_{\tau_k\le \la}|e_{\tau_k}(x)|^2\ls \la^n.
\end{equation} 
Since the eigenvalues of $H_V$ are all $\ge1$, by \eqref{rough} we get a crude long-time estimate
\begin{equation}\label{long}
	|e^{-tH_V}(x,y)|\ls e^{-t/2},\ \ t>1.
\end{equation}
For more information of the eigenfunction bounds on manifolds with boundary, see e.g. Smith-Sogge \cite{ssjams,ssacta}, Grieser \cite{Gri}, Sogge \cite{soggebd}, Koch-Smith-Tataru \cite{kst1,kst2},  Xu \cite{xu}.

\subsection{Main results} 
\begin{theorem}\label{thm1}
	Let $(M,g)$ be a smooth compact manifold of dimension $n\ge2$ with smooth boundary. Let $V\in L^1(M)$ with $V^-\in \mathcal{K}(M)$.  Then we have
	\begin{equation}\label{Vweyl}
		\# \{k: \, \tau_k\le \la\}=(2\pi)^{-n}\omega_n |M|\lambda^n+O(\lambda^{n-1})
	\end{equation}
	under either Dirichlet, Neumann or Robin boundary condition.
\end{theorem}

\begin{theorem}\label{thm2}
			Let $(M,g)$ be a smooth compact manifold of dimension $n\ge2$ with smooth boundary.  Suppose that the set of all periodic geodesic billiards has measure 0.  Let $V\in L^1(M)$ with $V^-\in \mathcal{K}(M)$.  Then we have  
	\begin{equation}\label{Vweylc}	\#\{k:\tau_k\le \la\}=(2\pi)^{-n}\omega_n|M|\la^n\mp\frac14(2\pi)^{1-n}\omega_{n-1}|\partial M|\la^{n-1}+o(\la^{n-1}),\end{equation}
	where the minus corresponds to the Dirichlet condition, and the plus to the Neumann condition.
\end{theorem}

 Recall that Huang-Sogge \cite{hs} proved Weyl laws for Schr\"odinger operators with critically singular potentials on compact manifolds without boundary. They also improved the Weyl remainder estimates under certain conditions, such as the condition \eqref{DGcond} by Duistermaat-Guillemin \cite{DG}. Frank-Sabin \cite{fs} obtained Weyl laws  for Schr\"odinger operators on 3-dimensional compact manifolds with or without boundary. They also proved sharp pointwise Weyl laws  on 3-dimensional compact manifolds without boundary. Inspired by these works, Huang-Zhang \cite{hz2021,hz2023} proved sharp pointwise Weyl laws for Schr\"odinger operators on compact manifolds of dimension $n\ge2$. 

In this paper, we simplify the perturbation argument in \cite{hs} and only use the Gaussian heat kernel bounds \eqref{heatV} for short times to control the remainder terms.  Indeed, the original argument  in \cite{hs} needs the Hadamard parametrix to calculate the kernels of $m(\sqrt{-\Delta_g})$, i.e. the functions of the Laplacian (see Sogge \cite[Section 4.3]{fio}). However, it is difficult to construct a precise parametrix for the wave kernel near the boundary, see e.g. Seeley \cite{seeley1, seeley2}, Ivrii \cite{ivrii}, Melrose-Taylor \cite{MTbook}, Smith-Sogge \cite{ssjams,ssacta} and H\"ormander \cite{hor4}. We can get around this difficulty, since our new argument does not involve the boundary, as long as the Gaussian heat kernel bounds \eqref{heatV} for short times are valid.

  \subsection{Paper structure} The paper is structured as follows. In Section 2, we introduce the  perturbation argument involving the wave equation to prove the theorems, and reduce to the perturbation estimates for short-interval spectral projection \eqref{localsum} and long-interval spectral projection \eqref{weylsum}. In Section 3 and 4, we prove these two estimates by frequency decomposition and the Gaussian heat kernel bounds \eqref{heatV} for short times.

   \subsection{Notations}Throughout this paper, $X\ls Y$ means $X\le CY$  for some positive constants $C$. The constant $C$ may depend on the potential $V$ and the manifold $(M,g)$, but it is independent of the parameter $\la$ and $\eps$. If $X\ls Y$ and $Y\ls X$, we denote $X\approx Y$.

   \section{Main argument}
   Let $0<\eps\le1$. Let $N^0(\la)=\#\{j:\la_j\le\la\}$ and $N_V(\la)=\#\{k:\tau_k\le \la\}$. We assume throughout that 
   \begin{equation}\label{ass}
   	N^0(\la)=c_0\la^n+c_1\la^{n-1}+O(\eps \la^{n-1}).
   \end{equation}
   Here $c_0,c_1$ are the coefficients in \eqref{weyllaw} or \eqref{weylconj}. 
   We shall prove that there exists $C_\eps>0$ such that
    \begin{equation}\label{goal}
   	N_V(\la)=c_0\la^n+c_1\la^{n-1}+O(\eps \la^{n-1}+C_\eps\la^{n-\frac32}).
   \end{equation}
    We shall fix $\eps=1$ to prove Theorem \ref{thm1}, and we choose $\eps$ arbitrarily small and $\la$ sufficiently large to prove Theorem \ref{thm2}. 
   
   We briefly review the perturbation argument in \cite{hs}.  We denote the indicator function of the interval $[-\la,\la]$ by $\ola(\tau)$. Then we can represent $N^0(\la)$ and $N_V(\la)$ as the trace of $\1_\la(P^0)$ and $\1_\la(P_V)$. Namely,
   \begin{equation}\label{N0}
   	N^0(\la)=\int_M\1_\la(P^0)(x,x)dx,\ \ \ N_V(\la)=\int_M\1_\la(P_V)(x,x)dx.
   \end{equation}
So to prove \eqref{goal} it suffices to estimate the trace of $\1_\la(P_V)-\1_\la(P^0)$. As is the custom (cf. \cite{hangzhou}),  we shall consider the $\eps$-dependent  approximation $\tilde\1_\la(\tau)$  (see \eqref{app}) and it suffices to prove the trace estimates
\begin{equation}\label{tr1}
	\Big|\int_M \big(\tilde\1_\la(P^0)-\1_\la(P^0)\big)(x,x)dx\Big|\ls \eps\la^{n-1},
\end{equation}\begin{equation}\label{tr2}
\Big|\int_M \big(\tilde\1_\la(P_V)-\1_\la(P_V)\big)(x,x)dx\Big|\ls \eps\la^{n-1}+C_\eps\la^{n-\frac32},
\end{equation}\begin{equation}\label{tr3}
\Big|\int_M \big(\tilde\1_\la(P_V)-\tilde \1_\la(P^0)\big)(x,x)dx\Big|\ls \eps\la^{n-1}+C_\eps\la^{n-\frac32}.
\end{equation}

Let $\rho\in C_0^\infty(\mathbb{R})$ be a fixed even real-valued function satisfying
\[\rho(t)=1\ \text{on}\ [-1/2,1/2]\ \text{and}\ \supp\rho\subset(-1,1).\] We define
\begin{equation}\label{app}
	\tilde\1_\la(\tau)=\frac1\pi\int_{\mathbb{R}}\rho(\eps t)\frac{\sin\la t}t\cos t\tau dt.
\end{equation}
 Since the Fourier transform of $\ola(\tau)$ is $2\frac{\sin\la t}t$, we have the  rapid decay property for $\tau\ge1$
  \begin{equation}\label{diff}
  	|\1_\la(\tau)-\tilde \1_\la(\tau)|\ls (1+\eps^{-1}|\la-\tau|)^{-N},\ \ \forall  N,
  \end{equation}
     \begin{equation}\label{diff2}
   	|\partial_\tau^j\tilde \1_\la(\tau)|\ls \eps^{-j} (1+\eps^{-1}|\la-\tau|)^{-N},\ \ \forall  N,\  \text{if}\ j=1,2,....
   \end{equation}

First, by \eqref{ass}, we have the pointwise estimate
\begin{equation}\label{local}
\#\{j:\la_j\in[\la,\la+\eps]\}=	\int_M\sum_{\la_j\in [\la,\la+\eps]}|e_j^0(x)|^2dx\ls \eps\la^{n-1}.
\end{equation}
Then we have the trace estimate 
\begin{align*}\Big|\int_M (\tilde\1_\la(P^0)-\1_\la(P^0))(x,x)dx\Big|&\ls \int_M\sum_{j}(1+\eps^{-1}|\la-\la_j|)^{-N}|e_j^0(x)|^2dx\ls \eps\la^{n-1}.
\end{align*}

Second, to handle the trace of $\tilde\1_\la(P_V)-\1_\la(P_V)$, it suffices to prove the short-interval estimate
\begin{equation}\label{localV}
\#\{k:\tau_k\in[\la,\la+\eps]\}=	\int_M\sum_{\tau_k\in [\la,\la+\eps]}|e_{\tau_k}(x)|^2dx\ls \eps\la^{n-1}+C_\eps\la^{n-\frac32}.
\end{equation}
Let $\chi\in C_0^\infty(\mathbb{R})$ be a nonnegative function with $\supp \chi\in (-1,1)$. Let $$\tilde\chi_\la(\tau)=\chi(\eps^{-1}(\la-\tau)).$$ For $\tau\ge1$ we have $\chi(\eps^{-1}(\la+\tau))=0$.
So for $\tau>0$ we obtain
\[\tilde \chi_\la(\tau)=\chi(\eps^{-1}(\la-\tau))+\chi(\eps^{-1}(\la+\tau))=\frac1\pi\int_{\mathbb{R}}\eps\hat\chi(\eps t)e^{it\la}\cos t\tau dt.\]
We  have for $\tau\ge1$,
\begin{equation}\label{chibd}
	|\partial_\tau^j\tilde \chi_\la(\tau)|\ls \eps^{-j}\1_{\{|\la-\tau|<\eps\}}(\tau),\ \ \forall N,\ \ \text{if}\ j=0,1,2,....
\end{equation}
To prove \eqref{localV}, it suffices to estimate the trace of $\tilde\chi_\la(P_V)$.  But the trace of $\tilde\chi_\la(P^0)$ satisfies the same bound as \eqref{local}. So we only need to handle the trace of their difference $\tilde\chi_\la(P_V)-\tilde\chi_\la(P^0)$. By Duhamel's principle and the spectral theorem, we can calculate the difference between the wave kernel and its perturbation (see \cite{hs}, \cite{hz2021}, \cite{hz2023})
\begin{align*}\cos tP_V(x,y)&-\cos tP^0(x,y)\\
	&=-\sum_{\la_j}\sum_{\tau_k}\int_M\int_0^t\frac{\sin(t-s)\la_j}{\la_j}\cos s\tau_k\  e_j^0(x)e_j^0(z)V(z)e_{\tau_k}(z)e_{\tau_k}(y)dzds\\
	&=\sum_{\la_j}\sum_{\tau_k}\int_M\frac{\cos t\la_j-\cos t\tau_k}{\la_j^2-\tau_k^2} e_j^0(x)e_j^0(z)V(z)e_{\tau_k}(z)e_{\tau_k}(y)dz.\end{align*}
Thus, it suffices to prove the short-interval estimate
     \begin{equation}\label{localsum}\Bigl| \, \sum_{\lambda_j}
   	\sum_{\tau_k}\int_M \int_M
   	\frac{\tilde \chi_\la(\la_j)-\tilde \chi_\la(\tau_k)}{\la_j^2-\tau_k^2} e_j^0(x)e_j^0(y)V(y)e_{\tau_k}(x)
   	e_{\tau_k}(y) \, dx dy\, \Bigr|\ls \eps\la^{n-1}+C_\eps \la^{n-\frac32}.\end{equation}
   	
 Third, to handle the trace of $\tilde\1_\la(P_V)-\tilde \1_\la(P^0)$, similarly it suffices to prove the long-interval estimate
  \begin{equation}\label{weylsum}\Bigl| \, \sum_{\lambda_j}
\sum_{\tau_k}\int_M \int_M
\frac{\ala(\la_j)-\ala(\tau_k)}{\la_j^2-\tau_k^2} e_j^0(x)e_j^0(y)V(y)e_{\tau_k}(x)
e_{\tau_k}(y) \, dx dy\, \Bigr|\ls \eps\la^{n-1}+C_\eps \la^{n-\frac32}.\end{equation}
In the following two sections, we shall prove \eqref{localsum} and \eqref{weylsum}.

\section{Proof of the short-interval estimate}
In this section, we shall prove the perturbation estimate for short-interval spectral projection. \begin{prop}Let $V\in L^1(M)$ with $V^-\in \mathcal{K}(M)$.  Then for any $\eps>0$ there exists $C_\eps>0$ such that
	\begin{equation}
		\#\{k:\tau_k\in[\la,\la+\eps]\}\ls \eps\la^{n-1}+C_\eps\la^{n-\frac32}.
	\end{equation}
\end{prop}
By the reduction in the Section 2, it suffices to prove \eqref{localsum}, namely
 \begin{equation}\label{localsum'}\Bigl| \, \sum_{\lambda_j}
	\sum_{\tau_k}\int_M \int_M
	\frac{\tilde \chi_\la(\la_j)-\tilde \chi_\la(\tau_k)}{\la_j^2-\tau_k^2} e_j^0(x)e_j^0(y)V(y)e_{\tau_k}(x)
	e_{\tau_k}(y) \, dx dy\, \Bigr|\ls \eps\la^{n-1}+C_\eps \la^{n-\frac32}.\end{equation}

 We shall split $V=V_0+V_1$ such that $V_0\in L^\infty(M)$ and $\|V_1\|_{L^1(M)}<\eps^2$. Let $$m(\la_j,\tau_k)=\frac{\tilde \chi_\la(\la_j)-\tilde \chi_\la(\tau_k)}{\la_j^2-\tau_k^2}.$$
 By the support property of $m(\la_j,\tau_k)$, we need to consider five cases.

\begin{enumerate}
	\item $|\tau_k-\la|\le \eps$, $|\la_j-\la|\le \eps$.
	\item $|\tau_k-\la|\le \eps$, $|\la_j-\la|\in(2^\ell,2^{\ell+1}]$, $\eps\le  2^\ell\le \la$.
	
	\item $|\la_j-\la|\le\eps$, $|\tau_k-\la|\in(2^\ell,2^{\ell+1}]$, $\eps\le 2^\ell\le \la$.
	\item $|\la_j-\la|\le \eps$, $\tau_k>2\la$.
	\item $|\tau_k-\la|\le \eps$, $\la_j>2\la$.
\end{enumerate}
We shall only use the Gaussian heat kernel bounds \eqref{heatV} for short times, and the argument is essentially symmetric in $\la_j$ and $\tau_k$. So the proofs of Case 3 and Case 5 are the the same as Case 2 and Case 4, respectively.

$ $

\noindent	\textbf{Case 1.} $|\tau_k-\la|\le \eps$, $|\la_j-\la|\le \eps$.

In this case, for $|s-\la|\le\eps$ we have
\[|m(\la_j,s)|+|\eps\partial_s m(\la_j,s)|\ls \eps^{-1}\la^{-1}.\]Then 
\begin{align*}&\sum_{|\la_j-\la|\le\eps}\sum_{|\tau_k-\la|\le \eps}\int_M\int_Mm(\la_j,\tau_k)e_j^0(x)e_j^0(y)V(y)e_{\tau_k}(y)e_{\tau_k}(x)dydx\\
	&= \sum_{|\la_j-\la|\le\eps}\sum_{|\tau_k-\la|\le \eps}\int_M\int_M\int_{\la-\eps}^{\la+\eps}\partial_sm(\la_j,s)\1_{[\la-\eps,\tau_k]}(s)e_j^0(x)e_j^0(y)V(y)e_{\tau_k}(y)e_{\tau_k}(x)dydxds\\
	&\ \ \ +\sum_{|\la_j-\la|\le\eps}\sum_{|\tau_k-\la|\le \eps}\int_M\int_Mm(\la_j,\la-\eps)e_j^0(x)e_j^0(y)V(y)e_{\tau_k}(y)e_{\tau_k}(x)dydx\\
	&=I_1+I_2
\end{align*}
We first handle $I_2$, and $I_1$ can be handled similarly. Since  we split $V=V_0+V_1$ such that $V_0\in L^\infty(M)$ and $\|V_1\|_{L^1(M)}<\eps^2$, we handle these two  parts separately. First, by H\"older inequality and the eigenfunction bound \eqref{rough} we have
\begin{align*}
	&\Big|\sum_{|\la_j-\la|\le\eps}\sum_{|\tau_k-\la|\le \eps}\int_M\int_Mm(\la_j,\la-\eps)e_j^0(x)e_j^0(y)V_1(y)e_{\tau_k}(y)e_{\tau_k}(x)dydx\Big|\\
	&\ls \|V_1\|_{L^1(M)}\cdot \eps^{-1}\la^{-1}\cdot \sup_y(\sum_{|\la_j-\la|\le \eps}|e_j^0(y)|^2)^\frac12(\sum_{|\tau_k-\la|\le \eps}|e_{\tau_k}(y)|^2)^\frac12\\
	&\ls \|V_1\|_{L^1(M)}\cdot \eps^{-1}\la^{-1}\cdot \la^{n/2}\cdot \la^{n/2}\\
	&\ls \eps\la^{n-1}.
\end{align*}
Similarly, by \eqref{local} and the eigenfunction bound \eqref{rough} we have
\begin{align*}
	&\Big|\sum_{|\la_j-\la|\le\eps}\sum_{|\tau_k-\la|\le \eps}\int_M\int_Mm(\la_j,\la-\eps)e_j^0(x)e_j^0(y)V_0(y)e_{\tau_k}(y)e_{\tau_k}(x)dydx\Big|\\
	&\ls \|V_0\|_{L^\infty(M)}\cdot \eps^{-1}\la^{-1}\cdot (\sum_{|\la_j-\la|\le \eps}1)^\frac12(\sum_{|\tau_k-\la|\le \eps}1)^\frac12\\
	&\ls \|V_0\|_{L^\infty(M)}\cdot \eps^{-1}\la^{-1}\cdot (\eps\la^{n-1})^{1/2}\cdot \la^{n/2}\\
	&\ls C_\eps \la^{n-\frac32}.
\end{align*} 
Combing these two parts, we get the desired bound.

$ $

\noindent \textbf{Case 2.} $|\tau_k-\la|\le \eps$, $|\la_j-\la|\in(2^\ell,2^{\ell+1}]$, $\eps\le  2^\ell\le \la$.

When $|\la_j-\la| \in(2^\ell,2^{\ell+1}]$, we have $m(\la_j,\tau_k)=\frac{-\tilde \chi_\la(\tau_k)}{\la_j^2-\tau_k^2}$, and for $|s-\la|\le \eps$ 
\begin{equation}\label{m2}|m(\la_j,s)|+|\eps\partial_s m(\la_j,s)|\ls  \la^{-1}2^{-\ell}.\end{equation}
We can use the same argument as Case 1 to handle 
\[m(\la_j,\tau_k)=m(\la_j,\la-\eps)+\int_{\la-\eps}^{\la+\eps}\partial_sm(\la_j,s)\1_{[\la-\eps,\tau_k]}(s)ds.\]
As before, we just need to handle the first term, and the second term is similar by \eqref{m2}. Since   $V=V_0+V_1$ with $V_0\in L^\infty(M)$ and $\|V_1\|_{L^1(M)}<\eps^2$, we shall handle these two  parts separately. First, by H\"older inequality and the eigenfunction bound \eqref{rough} we have
\begin{align*}
	&\Big|\sum_{|\la_j-\la|\in(2^\ell,2^{\ell+1}]}\sum_{|\tau_k-\la|\le \eps}\int_M\int_Mm(\la_j,\la-\eps)e_j^0(x)e_j^0(y)V_1(y)e_{\tau_k}(y)e_{\tau_k}(x)dydx\Big|\\
	&\ls \|V_1\|_{L^1(M)}\cdot 2^{-\ell}\la^{-1}\cdot \sup_y(\sum_{|\la_j-\la|\in(2^\ell,2^{\ell+1}]}|e_j^0(y)|^2)^\frac12(\sum_{|\tau_k-\la|\le \eps}|e_{\tau_k}(y)|^2)^\frac12\\
	&\ls \|V_1\|_{L^1(M)}\cdot 2^{-\ell}\la^{-1}\cdot \la^{n/2}\cdot \la^{n/2}\\
	&\ls \eps^22^{-\ell}\la^{n-1}.
\end{align*}
Similarly, by \eqref{local} and the eigenfunction bound \eqref{rough} we have
\begin{align*}
	&\Big|\sum_{|\la_j-\la|\in(2^\ell,2^{\ell+1}]}\sum_{|\tau_k-\la|\le \eps}\int_M\int_Mm(\la_j,\la-\eps)e_j^0(x)e_j^0(y)V_0(y)e_{\tau_k}(y)e_{\tau_k}(x)dydx\Big|\\
	&\ls \|V_0\|_{L^\infty(M)}\cdot 2^{-\ell}\la^{-1}\cdot (\sum_{|\la_j-\la|\in(2^\ell,2^{\ell+1}]}1)^\frac12(\sum_{|\tau_k-\la|\le \eps}1)^\frac12\\
	&\ls \|V_0\|_{L^\infty(M)}\cdot 2^{-\ell}\la^{-1}\cdot (2^\ell\la^{n-1})^{1/2}\cdot \la^{n/2}\\
	&\ls C_\eps 2^{-\ell/2}\la^{n-\frac32}.
\end{align*} 
Summing over $\ell\in \mathbb{Z}:\ \eps\le 2^\ell\le \la$, we get the desired bound.

$ $

\noindent \textbf{Case 3.} $|\la_j-\la|\le \eps$, $|\tau_k-\la|\in(2^\ell,2^{\ell+1}]$, $\eps\le  2^\ell\le \la$.

This case is essentially the same as Case 2, since we only use the Gaussian heat kernel bound and the proof still works if we interchange $\la_j$ and $\tau_k$.

$ $

\noindent \textbf{Case 4.} $|\la_j-\la|\le \eps$, $\tau_k>2\la$.

In this case, $m(\la_j,\tau_k)=\frac{\tilde \chi_\la(\la_j)}{\la_j^2-\tau_k^2}$. We expand
\begin{align*}
	\frac1{\tau_k^2-\lambda_j^2}=\tau_k^{-2}+\tau_k^{-2}(\la_j/\tau_k)^2+...+\tau_k^{-2}(\la_j/\tau_k)^{2N-2}+\frac{(\la_j/\tau_k)^{2N}}{\tau_k^2-\la_j^2}.
\end{align*}
We will fix $N=2n$ later. For $\ell=0,1,...,N-1$, when $n-4-4\ell<0$ we use \eqref{rough} to get
\begin{align*}
	&\Big|\sum_{|\la_j-\la|\le \eps}\sum_{\tau_k>2\la}\int_M\int_M\tilde\chi_\la(\la_j)\la_j^{2\ell}e_j^0(x)e_j^0(y)V(y)\tau_k^{-2-2\ell}e_{\tau_k}(y)e_{\tau_k}(x)dxdy\Big|\\
	&\ls \|V\|_{L^1(M)}\cdot \la^{2\ell}\sup_y(\sum_{|\la_j-\la|\le \eps}|e_j^0(y)|^2)^{1/2}(\sum_{\tau_k>2\la}\tau_k^{-4-4\ell}|e_{\tau_k}(y)|^2)^{1/2}\\
	&\ls \|V\|_{L^1(M)}\cdot \la^{2\ell}\cdot\la^{n/2}\cdot \la^{\frac n2-2-2\ell}\\
	&\ls \la^{n-2}.
\end{align*}
When $n-4-4\ell\ge0$, we split the sum over $\tau_k>2\la$ into the difference between the complete sum and the partial sum $\tau_k\le 2\la$. We first handle the partial sum by \eqref{rough}. Then we have
\begin{align*}
	&\Big|\sum_{|\la_j-\la|\le \eps}\sum_{\tau_k\le 2\la}\int_M\int_M\tilde\chi_\la(\la_j)\la_j^{2\ell}e_j^0(x)e_j^0(y)V(y)\tau_k^{-2-2\ell}e_{\tau_k}(y)e_{\tau_k}(x)dxdy\Big|\\
	&\ls \|V\|_{L^1(M)}\cdot \la^{2\ell}\sup_y(\sum_{|\la_j-\la|\le \eps}|e_j^0(y)|^2)^{1/2}(\sum_{\tau_k\le2\la}\tau_k^{-4-4\ell}|e_{\tau_k}(y)|^2)^{1/2}\\
	&\ls \|V\|_{L^1(M)}\cdot \la^{2\ell}\cdot\la^{n/2}\cdot \la^{\frac n2-2-2\ell}(\log\la)^\frac12\\
	&\ls \la^{n-2}(\log\la)^\frac12.
\end{align*}
The factor $(\log\la)^\frac12$ only appears when $n-4-4\ell=0$. We shall handle the complete sum by the kernel estimates
\begin{equation}
	|\tilde\chi_\la(P^0)(P^0)^{2\ell}(x,y)|\ls \la^{n+2\ell},
\end{equation}
\begin{equation}
	\|\tilde\chi_\la(P^0)(P^0)^{2\ell}(\cdot,y)\|_{L^2(M)}\ls \la^{\frac n2+2\ell},
\end{equation}
and when $n-4-4\ell\ge0$, 
\begin{equation}\label{rieszk}
	|H_V^{-1-\ell}(x,y)|\ls d_g(x,y)^{-n+2+2\ell}.
\end{equation}
These follow from the heat kernel bounds \eqref{heatV} and \eqref{long}, and the relation
\[H_V^{-1-\ell}(x,y)=\frac1{\ell!}\int_0^\infty t^\ell e^{-tH_V}(x,y)dt.\]
Then by H\"older inequality we have
\begin{align*}
	&\Big|\sum_{|\la_j-\la|\le \eps}\sum_{\tau_k}\int_M\int_M\tilde\chi_\la(\la_j)\la_j^{2\ell}e_j^0(x)e_j^0(y)V(y)\tau_k^{-2-2\ell}e_{\tau_k}(y)e_{\tau_k}(x)dxdy\Big|\\
	&\ls \|V\|_{L^1(M)}\cdot \sup_y\int_M|\tilde\chi_\la(P^0)(P^0)^{2\ell}(x,y)H_V^{-1-\ell}(x,y)|dx\\
	&\ls  \sup_y \int_{d_g(x,y)\le\la^{-1}}\la^{n+2\ell}d_g(x,y)^{-n+2+2\ell}dx+\la^{\frac n2+2\ell}\sup_y\Big(\int_{d_g(x,y)>\la^{-1}}d_g(x,y)^{-2n+4+4\ell}dx\Big)^{1/2}\\
	&\ls \la^{n-2}(\log\la)^\frac12.
\end{align*}
The factor $(\log\la)^\frac12$ only appears when $n-4-4\ell=0$.

Now we handle the last term in the expansion. Let $$m_N(\la_j,s)=\tilde\chi_\la(\la_j)\frac{\la_j^{2N}}{1-s^2\la_j^2}.$$ Then for $s\in [0,(2\la)^{-1}]$ we have
\begin{equation}\label{mN0}|m_N(\la_j,s)|+|\la^{-1}\partial_s m_N(\la_j,s)|\ls \la^{2N}.\end{equation}
We can use the same argument as Case 1 to handle 
\[m_N(\la_j,\tau_k^{-1})=m_N(\la_j,0)+\int_0^{(2\la)^{-1}}\partial_s m_N(\la_j,s)\1_{[0,1/\tau_k]}(s)ds.\]
As before, we just need to handle the first term, and the second term is similar by \eqref{mN0}. For $N=2n$, by using the eigenfunction bound \eqref{rough} we have 
\begin{align*}
	&\Big|\sum_{|\la_j-\la|\le \eps}\sum_{\tau_k>2\la}\int_M\int_M\tilde\chi_\la(\la_j)\la_j^{2N}e_j^0(x)e_j^0(y)V(y)\tau_k^{-2-2N}e_{\tau_k}(y)e_{\tau_k}(x)dxdy\Big|\\
	&\ls \|V\|_{L^1(M)}\cdot \la^{2N}\sup_y(\sum_{|\la_j-\la|\le \eps}|e_j^0(y)|^2)^{1/2}(\sum_{\tau_k>2\la}\tau_k^{-4-4N}|e_{\tau_k}(y)|^2)^{1/2}\\
	&\ls \|V\|_{L^1(M)}\cdot \la^{2N}\cdot\la^{n/2}\cdot \la^{\frac n2-2-2N}\\
	&\ls \la^{n-2}.
\end{align*}
$ $

\noindent \textbf{Case 5.} $|\tau_k-\la|\le \eps$, $\la_j>2\la$.

This case is essentially the same as Case 4, since we only use the Gaussian heat kernel bound and the proof still works if we interchange $\la_j$ and $\tau_k$.
\section{Proof of the long-interval estimate}
In this section, we prove the perturbation estimate \eqref{weylsum} for the long-interval spectral projection . The argument is essentially similar to the proof of the short-interval estimate \eqref{localsum}.

\begin{prop}Let $V\in L^1(M)$ with $V^-\in \mathcal{K}(M)$.  Then for any $\eps>0$ there exists $C_\eps>0$ such that
 \begin{equation}\label{weylsum'}\Bigl| \, \sum_{\lambda_j}
	\sum_{\tau_k}\int_M \int_M
	\frac{\ala(\la_j)-\ala(\tau_k)}{\la_j^2-\tau_k^2} e_j^0(x)e_j^0(y)V(y)e_{\tau_k}(x)
	e_{\tau_k}(y) \, dx dy\, \Bigr|\ls \eps\la^{n-1}+C_\eps \la^{n-\frac32}.\end{equation}
\end{prop}
 If $\mu$ is the frequency, then we denote Low=``$\mu<\lambda/2$'', Med=``$\lambda/2\le \mu\le 10\la$'', High=``$\mu>10\lambda$''. We shall split the sum into the following five cases:
\begin{enumerate}
	\item Low+Low: $\la_j<\lambda/2$ and $\tau_k<\lambda/2$
	\item MedLow+Med: $\la_j\le 10\la$ and $\tau_k\in[\la/2,10\la]$
	\item Med+Low: $\la_j\in[\la/2,10\la]$ and $\tau_k<\la/2$
	\item All+High: all $\la_j$, and $\tau_k>10\la$
	\item High+MedLow: $\la_j>10\la$ and $\tau_k\le 10\la$.
\end{enumerate}
\subsection{Low+Low}
\begin{prop}\label{lowlow}
	Let $V\in L^1(M)$ with $V^-\in \mathcal{K}(M)$.   Then
	\begin{align}\label{llow}
		\Bigl| \, \sum_{\la_j<\lambda/2}
		\sum_{ \tau_k<\lambda/2}\int_M \int_M
		\frac{\ala(\la_j)-\ala(\tau_k)}{\la_j^2-\tau_k^2} e_j^0(x)e_j^0(y)V(y)e_{\tau_k}(x)
		e_{\tau_k}(y) \, dx dy\, \Bigr|
		\ls \la^{-\sigma},\ \forall \sigma.
	\end{align}
\end{prop}

This case simply follows from  the rapid decay property \eqref{chibd} and the eigenfunction bound \eqref{rough}. Indeed, by the mean value theorem and  \eqref{chibd} we have
\[\Big|\frac{\ala(\la_j)-\ala(\tau_k)}{\la_j^2-\tau_k^2}\Big|\ls \eps^{-1}(1+\eps^{-1}\la)^{-N},\ \forall N.\]
Thus, by using \eqref{rough} we obtain
\begin{align*}
	\text{L.H.S. of}\ \eqref{llow}\ls \eps^{-1}(1+\eps^{-1}\la)^{-N}\cdot \la^n\cdot \la^n,\ \ \forall N.
\end{align*}
This implies \eqref{llow}.

\subsection{MedLow+Med, Med+Low}
\begin{prop}\label{medmed}
	Let $V\in L^1(M)$ with $V^-\in \mathcal{K}(M)$.   Then
	\begin{align}\label{ml1}
		\Bigl| \, \sum_{\lambda_j\le10\lambda}
		\sum_{\tau_k\in [\lambda/2,10\lambda]}\int_M \int_M
		\frac{\ala(\la_j)-\ala(\tau_k)}{\la_j^2-\tau_k^2} e_j^0(x)e_j^0(y)V(y)e_{\tau_k}(x)
		e_{\tau_k}(y) \, dx dy\, \Bigr|\ls \eps
		\la^{n-1}+C_\eps \la^{n-\frac32},
	\end{align}
	\begin{align}\label{ml2}
		\Bigl| \, \sum_{\tau_k<\lambda/2}
		\sum_{\la_j\in [\lambda/2,10\lambda]}\int_M \int_M
		\frac{\ala(\la_j)-\ala(\tau_k)}{\la_j^2-\tau_k^2} e_j^0(x)e_j^0(y)V(y)e_{\tau_k}(x)
		e_{\tau_k}(y) \, dx dy\, \Bigr|\ls
		\eps\la^{n-1}+C_\eps \la^{n-\frac32}.
	\end{align}
\end{prop}

We just need to handle \eqref{ml1}, and the proof of \eqref{ml2} is similar, since we only use the eigenfunction bound \eqref{rough} and the proof of \eqref{ml1} still works if we interchange $\la_j$ and $\tau_k$.

 We fix a Littlewood-Paley bump function
$\beta\in C^\infty_0((1/2,2))$ satisfying
$$
\sum_{\ell=-\infty}^\infty \beta(2^{-\ell} s)=1, \quad s>0.$$
Let $\ell_0\le 0$ be the largest integer such that $2^{\ell_0}\le \eps$, and let
$$\beta_0(s)=\sum_{\ell \le \ell_0} \beta(2^{-\ell}|s|)
\in C^\infty_0((-2,2)),$$
and
$$\tilde \beta(s)=s^{-1}\beta(|s|)\in C^\infty_0
\bigl(\{|s|\in (1/2,2)\}\bigr).$$
To prove \eqref{ml1}, we write for $\la_j\le 10\la$ and $\la/2\le \tau\le 10\la$
\begin{align*}
m(\la_j,\tau)&= \frac{\ala(\la_j)-\ala(\tau)}{\la_j^2-\tau^2}
\\
&=  \frac{\ala(\la_j)-\ala(\tau)}{\la_j-\tau}\frac{\beta_0(\la_j-\tau)}{\la_j+\tau} +\sum_{\eps< 2^\ell \ls \la}
\frac{2^{-\ell}\tilde\beta(2^{-\ell}(\la_j-\tau))}
{\la_j+\tau}
\, (\ala(\la_j)-\ala(\tau)) .
\end{align*}
We let
$$m_0(\la_j,\tau)=  \frac{\ala(\la_j)-\ala(\tau)}{\la_j-\tau}\frac{\beta_0(\la_j-\tau)}{\la_j+\tau}$$
and if  $2^\ell \ls \la$, let
$$R_{\ell}(\la_j,\tau)= 
\frac{2^{-\ell}\tilde\beta(2^{-\ell}(\la_j-\tau))}
{\la_j+\tau}$$
\begin{align*}
	m^-_{\ell}(\la_j,\tau)&=
	\frac{2^{-\ell}\tilde \beta(2^{-\ell}(\la_j-\tau))}{\la_j+\tau} \bigl(\ala(\la_j)-1\bigr)
	\\
	m^+_{\ell}(\la_j,\tau)&=
	\frac{2^{-\ell}\tilde \beta(2^{-\ell}(\la_j-\tau))}{\la_j+\tau} \ala(\la_j)
\end{align*}
So when $\tau\in[\la/2,\la]$, we can write
\[
m(\la_j,\tau)=m_0(\la_j,\tau)
+ \sum_{\eps<2^\ell \ls \la}\Big(m^-_{\ell}(\la_j,\tau)+R_{\ell}(\la_j,\tau)\bigl(1-\ala(\tau)\bigr)\Big)\]
and when $\tau\in(\la,10\la]$, we can write
\[
m(\la_j,\tau)=m_0(\la_j,\tau)
+ \sum_{\eps<2^\ell \ls \la}\Big(m^+_{\ell}(\la_j,\tau)-R_{\ell}(\la_j,\tau)\ala(\tau)\Big).\]
For  $\ell\ge\ell_0$ and $\nu=0,1,2,...$, let
\[I_{\ell,\nu}^-=(\la-(\nu+1)2^\ell,\la-\nu2^\ell]\ \ \text{and}\ \ I_{\ell,\nu}^+=(\la+\nu2^\ell,\la+(\nu+1)2^\ell].\]
 By the rapid decay property \eqref{diff} and \eqref{diff2}, we have
\begin{lemma}\label{symb}
	If $\ell\in \mathbb{Z}:\eps<2^\ell\ls \la$,
and $\nu=0,1,2,\dots$, 
\begin{equation}
	|m^{\pm}_{\ell}(\la_j,\tau)|+|2^\ell\partial_\tau m^\pm_{\ell}(\la_j,\tau)|\ls 2^{-\ell}\la^{-1}(1+\nu)^{-N},\ \ \tau\in I_{\ell,\nu}^\pm \cap [\la/2,10\la]
\end{equation}
\begin{equation}	|m_{0}(\la_j,\tau)|+|\eps\partial_\tau m_0(\la_j,\tau)|\ls \eps^{-1}\la^{-1}(1+\nu)^{-N},\ \ \tau\in I_{\ell_0,\nu}^\pm \cap [\la/2,10\la]
\end{equation}
\begin{equation}
	|R_\ell(\la_j,\tau)|+|\partial_\tau R_\ell(\la_j,\tau)|\ls 2^{-\ell}\la^{-1},\ \ \tau\in[\la/2,10\la].
\end{equation}
\end{lemma}
If we denote the left endpoint of the interval $I_{\ell,\nu}^\pm$ by $\tau_{\ell,\nu}^\pm$, then
\[m_\ell^\pm(\la_j,\tau)=m_\ell^\pm(\la_j,\tau_{\ell,\nu}^\pm)+\int_{I_{\ell,\nu}^\pm}\partial_s m_\ell^\pm(\la_j,s)\1_{[\tau_{\ell,\nu}^\pm,\tau]}(s)ds.\]
As before, we just need to handle the first term, and the second term is similar by Lemma \ref{symb}. We split $V=V_0+V_1$ such that $V_0\in L^\infty(M)$ and $\|V_1\|_{L^1(M)}<\eps^2$. Then by H\"older inequality and \eqref{rough} we have 
\begin{align*}
&\Big|\sum_{\la_j\le10\la}\sum_{\tau_k\in I_{\ell,\nu}^\pm}\int_M \int_M
m_\ell^\pm(\la_j,\tau_{\ell,\nu}^\pm) e_j^0(x)e_j^0(y)V_1(y)e_{\tau_k}(x)
e_{\tau_k}(y)dx dy\Big|\\
&\ls \|V_1\|_{L^1(M)}\cdot 2^{-\ell}\la^{-1}(1+\nu)^{-N}\cdot \sup_y(\sum_{|\la_j-\la|\ls (1+\nu) 2^\ell}|e_j^0(y)|^2)^{1/2}(\sum_{|\tau_k-\la|\ls (1+\nu) 2^\ell}|e_{\tau_k}(y)|^2)^{1/2}\\
&\ls \|V_1\|_{L^1(M)}\cdot 2^{-\ell}\la^{-1}(1+\nu)^{-N}\cdot \la^{n/2}\cdot \la^{n/2}\\
&\ls \eps^22^{-\ell}\la^{n-1}(1+\nu)^{-N}.
\end{align*}
Similarly, by \eqref{local} we have
\begin{align*}
	&\Big|\sum_{\la_j\le10\la}\sum_{\tau_k\in I_{\ell,\nu}^\pm}\int_M \int_M
	m_\ell^\pm(\la_j,\tau_{\ell,\nu}^\pm) e_j^0(x)e_j^0(y)V_0(y)e_{\tau_k}(x)
	e_{\tau_k}(y)dx dy\Big|\\
	&\ls \|V_0\|_{L^\infty(M)}\cdot 2^{-\ell}\la^{-1}(1+\nu)^{-2N}\cdot (\sum_{|\la_j-\la|\ls (1+\nu) 2^\ell}1)^{1/2}(\sum_{|\tau_k-\la|\ls (1+\nu) 2^\ell}1)^{1/2}\\
	&\ls \|V_0\|_{L^\infty(M)}\cdot 2^{-\ell}\la^{-1}(1+\nu)^{-2N}\cdot ( (1+\nu)2^\ell\la^{n-1})^{1/2}\cdot \la^{n/2}\\
	&\ls C_\eps2^{-\ell/2}\la^{n-\frac32}(1+\nu)^{-N}.
\end{align*}
By the same method, we can also obtain
\begin{align*}
	&\Big|\sum_{\la_j\le10\la}\sum_{\tau_k\in I_{\ell,\nu}^-}\int_M \int_M
R_\ell(\la_j,\tau_k)(1-\ala(\tau_k)) e_j^0(x)e_j^0(y)V(y)e_{\tau_k}(x)
	e_{\tau_k}(y)dx dy\Big|\\
	&\ls \eps^22^{-\ell}\la^{n-1}(1+\nu)^{-N}+C_\eps2^{-\ell/2}\la^{n-\frac32}(1+\nu)^{-N},
\end{align*}
\begin{align*}
	&\Big|\sum_{\la_j\le10\la}\sum_{\tau_k\in I_{\ell,\nu}^+}\int_M \int_M
	R_\ell(\la_j,\tau_k)\ala(\tau_k) e_j^0(x)e_j^0(y)V(y)e_{\tau_k}(x)
	e_{\tau_k}(y)dx dy\Big|\\
	&\ls \eps^22^{-\ell}\la^{n-1}(1+\nu)^{-N}+C_\eps2^{-\ell/2}\la^{n-\frac32}(1+\nu)^{-N},
\end{align*}
and
\begin{align*}
	&\Big|\sum_{\la_j\le10\la}\sum_{\tau_k\in I_{\ell_0,\nu}^\pm}\int_M \int_M
	m_0(\la_j,\tau_k) e_j^0(x)e_j^0(y)V(y)e_{\tau_k}(x)
	e_{\tau_k}(y)dx dy\Big|\\
	&\ls \eps\la^{n-1}(1+\nu)^{-N}+C_\eps\la^{n-\frac32}(1+\nu)^{-N}.
\end{align*}
Summing over $\ell\ge\ell_0$ and $\nu\ge0$ we get the desired bound.

\subsection{All+High, High+MedLow}
\begin{prop}\label{highfull}
Let $V\in L^1(M)$ with $V^-\in \mathcal{K}(M)$.   Then
\begin{equation}\label{ah1}
\Bigl| \, \sum_{\lambda_j}
\sum_{\tau_k >10\lambda}\int_M \int_M
\frac{\ala(\la_j)-\ala(\tau_k)}{\la_j^2-\tau_k^2} e_j^0(x)e_j^0(y)V(y)e_{\tau_k}(x)
e_{\tau_k}(y) \, dx dy\, \Bigr|
\ls \la^{n-2}(\log \la)^{1/2},
\end{equation}
\begin{equation}\label{ah2}
\Bigl| \, \sum_{\tau_k\le 10\la}
\sum_{\lambda_j >10\lambda}\int_M \int_M
\frac{\ala(\la_j)-\ala(\tau_k)}{\la_j^2-\tau_k^2} e_j^0(x)e_j^0(y)V(y)e_{\tau_k}(x)
e_{\tau_k}(y) \, dx dy\, \Bigr|
\ls
\la^{n-2}(\log \la)^{1/2}.
\end{equation}
\end{prop}

We just need to handle \eqref{ah1}, and the proof of \eqref{ah2} is similar, since we only use the Gaussian heat kernel bounds \eqref{heatV} for short times and the proof of \eqref{ah1} can still work if we interchange $\la_j$ and $\tau_k$. 

We first note that by the mean value theorem, the rapid decay property \eqref{diff2}, and the eigenfunction bound \eqref{rough} we have
\begin{align*}
\Bigl| \, \sum_{\la_j\in[\tau_k/2,2\tau_k]}
\sum_{\tau_k >10\lambda}\iint
\frac{\ala(\la_j)-\ala(\tau_k)}{\la_j^2-\tau_k^2}
e_j^0(x)e_j^0(y)V(y)e_{\tau_k}(x)
e_{\tau_k}(y) \, dx dy\, \Bigr|
\lesssim \la^{-\sigma} , \quad \forall \, \sigma.
\end{align*}
So now we only need to deal with two cases $\lambda_j<\tau_k/2$ and $\lambda_j>2\tau_k$.

\noindent \textbf{Case 1:} If $\lambda_j<\tau_k/2$, then by the rapid decay property \eqref{diff} we also have
\begin{align*}
\Bigl| \, \sum_{\la_j<\tau_k/2}
\sum_{\tau_k >10\lambda}\iint
\frac{\ala(\tau_k)}{\la_j^2-\tau_k^2}
e_j^0(x)e_j^0(y)V(y)e_{\tau_k}(x)
e_{\tau_k}(y) \, dx dy\, \Bigr|
\ls \la^{-\sigma}, \quad \forall \, \sigma.
\end{align*}
So we only need to prove
\begin{equation}\label{case1}
\Bigl| \, \sum_{\la_j<\tau_k/2}
\sum_{\tau_k >10\lambda}\iint
\frac{\ala(\lambda_j)}{\la_j^2-\tau_k^2}
e_j^0(x)e_j^0(y)V(y)e_{\tau_k}(x)
e_{\tau_k}(y) \, dx dy\, \Bigr|
\lesssim \la^{n-2} (\log\lambda)^\frac12.
\end{equation}
We expand
$$
\frac{1}{\tau_k^2-\la_j^2}=\tau_k^{-2}
+\tau_k^{-2}\bigl(\la_j/\tau_k\bigr)^2+
\cdots + \tau_k^{-2}\bigl(\la_j/\tau_k\bigr)^{2N-2}
\\
+(\la_j/\tau_k)^{2N}\frac{1}{\tau_k^2-\la_j^2}.
$$
where we will choose $N=2n$ later. It suffices to prove for $\ell=0,\dots,N-1,$
\begin{align}\label{3.23}
\Bigl| \iint \sum_{\tau_k>10\la}\tau_k^{-2-2\ell} e_{\tau_k}(x)e_{\tau_k}(y)\sum_{\lambda_j<\tau_k/2}
\la_j^{2\ell}
\ala(\la_j)e_j^0(x)e_j^0(y)
 \,
V(y) \, dxdy\Bigr|
\lesssim 
\la^{n-2} (\log\la)^{1/2},
\end{align}
as well as
\begin{align}\label{3.24}
\Bigl|\sum_{\tau_k>10\la}\sum_{\lambda_j<\tau_k/2}\iint
\frac{\la_j^{2N}}{\la_j^2-\tau_k^2}
\ala(\la_j) e^0_j(x)e^0_j(y)V(y)\tau_k^{-2N}
e_{\tau_k}(x)e_{\tau_k}(y) \, dxdy\Bigr|
\lesssim \la^{n-2}  .
\end{align}
First, we note that \eqref{3.23} is a consequence of
\begin{align}\label{3.23'}
\Bigl| \iint \bigl((P^0)^{2\ell}\ala(P^0)\bigr)(x,y)
\sum_{\tau_k>10\la}\tau_k^{-2-2\ell}
\, e_{\tau_k}(x)e_{\tau_k}(y) \,
V(y) \, dxdy\Bigr|
\lesssim 
\la^{n-2} (\log\la)^{1/2},
\end{align}
since if $\la_j\ge\tau_k/2>5\lambda$ then
\begin{align*}
\bigl|\ala(\la_j)\sum_{10\la<\tau_k\le 2\la_j}
\tau_k^{-2-2\ell}e_{\tau_k}(x)e_{\tau_k}(y)|
&\lesssim \la_j^{-\sigma}\sum_{10\la<\tau_k\le 2 \la_j}
\tau_k^{-\sigma}|\tau_k^{-2-2\ell}e_{\tau_k}(x)e_{\tau_k}(y)|
\\
&\lesssim \la_j^{n-2\ell-2\sigma},\ \ \ \forall \sigma
\end{align*}
which yields
\begin{align*}
\Bigl| \iint \sum_{\tau_k>10\la}\tau_k^{-2-2\ell} e_{\tau_k}(x)e_{\tau_k}(y)\sum_{\lambda_j\ge\tau_k/2}
\la_j^{2\ell}
\ala(\la_j)e_j^0(x)e_j^0(y)
 \,
V(y) \, dxdy\Bigr|\lesssim \la^{-\sigma},\ \ \forall \sigma.
\end{align*}
When $n-4-4\ell<0$, then by the eigenfunction bound \eqref{rough} we have
\begin{align*}
\text{L.H.S. of}\ \eqref{3.23'}&\ls \|V\|_{L^1(M)}\cdot \sup_y (\sum_{\la_j}\la_j^{4\ell}|\ala(\la_j)|^2|e_j^0(y)|^2)^{1/2}(\sum_{\tau_k>10\la}\tau_k^{-4-4\ell}|e_{\tau_k}(y)|^2)^{1/2}\\
&\lesssim \|V\|_{L^1} \cdot \la^{\frac{n}2+2\ell} \cdot \la^{\frac{n}2-2-2\ell}
\\
&\ls \la^{n-2}
\end{align*}
When $n-4-4\ell\ge0$, we split the sum over $\tau_k>10\la$ into the difference between the complete sum and the partial sum $\tau_k\le 10\la$. We first handle the partial sum by \eqref{rough}.  Then by the eigenfunction bound \eqref{rough} we have 
\begin{align*}
	&\Big|\iint \bigl((P^0)^{2\ell}\ala(P^0)\bigr)(x,y)\sum_{\tau_k\le 10\la}\tau_k^{-2-2\ell}e_{\tau_k}(y)e_{\tau_k}(x)V(y)dxdy\Big|\\
	&\ls \|V\|_{L^1(M)}\cdot \sup_y(\sum_{\la_j}\la^{4\ell}|\ala(\la_j)|^2|e_j^0(y)|^2)^{1/2}(\sum_{\tau_k\le10\la}\tau_k^{-4-4\ell}|e_{\tau_k}(y)|^2)^{1/2}\\
	&\ls \|V\|_{L^1(M)}\cdot\la^{\frac n2+2\ell}\cdot \la^{\frac n2-2-2\ell}(\log\la)^\frac12\\
	&\ls \la^{n-2}(\log\la)^\frac12.
\end{align*}
The factor $(\log\la)^\frac12$ only appears when $n-4-4\ell=0$.

Next, we handle the complete sum by the kernel estimates
\begin{equation}
	|((P^0)^{2\ell}\ala(P^0))(x,y)|\ls \la^{n+2\ell},
\end{equation}
\begin{equation}
	\|((P^0)^{2\ell}\ala(P^0))(\cdot,y)\|_{L^2(M)}\ls \la^{\frac n2+2\ell},
\end{equation}
and when $n-4-4\ell\ge0$, we recall \eqref{rieszk}, namely
\begin{equation}
	|H_V^{-1-\ell}(x,y)|\ls d_g(x,y)^{-n+2+2\ell}.
\end{equation}
Then by H\"older inequality we have
\begin{align*}
	&\Big|\iint((P^0)^{2\ell}\ala(P^0))(x,y)\sum_{\tau_k}\tau_k^{-2-2\ell}e_{\tau_k}(y)e_{\tau_k}(x)V(y)dxdy\Big|\\
	&\ls \|V\|_{L^1(M)}\cdot \sup_y\int_M|((P^0)^{2\ell}\ala(P^0))(x,y)H_V^{-1-\ell}(x,y)|dx\\
	&\ls  \sup_y \int_{d_g(x,y)\le\la^{-1}}\la^{n+2\ell}d_g(x,y)^{-n+2+2\ell}dx+\la^{\frac n2+2\ell}\sup_y\Big(\int_{d_g(x,y)>\la^{-1}}d_g(x,y)^{-2n+4+4\ell}dx\Big)^{1/2}\\
	&\ls \la^{n-2}(\log\la)^\frac12.
\end{align*}
The factor $(\log\la)^\frac12$ only appears when $n-4-4\ell=0$.

To prove \eqref{3.24}, we first note that if $N=2n$,
\begin{align*}
\Bigl|\sum_{\lambda_j<\tau_k/2}\sum_{\tau_k\ge \la^2}\iint
\frac{\la_j^{2N}}{\la_j^2-\tau_k^2}\ala(\la_j)
 e^0_j(x)e^0_j(y)V(y)\tau_k^{-2N}
e_{\tau_k}(x)e_{\tau_k}(y) \, dxdy\Bigr|\ls \la^{-N},
\end{align*}
since
\begin{align*}
\int_M\bigl| \, \sum_{\la_j<\tau_k/2}
\frac{\la_j^{2N}}{\la_j^2-\tau_k^2}
\ala(\la_j) e^0_j(x)e^0_j(y) \, \bigr| \, dx
&\lesssim \bigl\| \sum_{\la_j<\tau_k/2}
\frac{\la_j^{2N}}{\la_j^2-\tau_k^2}
\ala(\la_j) e^0_j(\, \cdot \, )e^0_j(y)\bigr\|_{L^2(M)}
\\
&\lesssim \bigl\| (P^0)^{2N} \ala(P^0)(\, \cdot ,y)\bigr\|_{L^2(M)} \lesssim \la^{\frac{n}2+2N}
\end{align*}
and $$\sum_{\tau_k\ge \la^2} \tau_k^{-2N}|e_{\tau_k}(x)e_{\tau_k}(y)|\ls \la^{-4N}\la^{2n}.$$
So we only need to handle the sum with $10\lambda<\tau_k<\la^2$. If $2\la<\la_j<\tau_k/2$, then
$$\frac{\la_j^{2N}}{\la^2_j-\tau_k^2}\ala(\la_j)=O(\tau_k^{-\sigma}),
\, \, \, \text{when} \, \, 10\la\le \tau_k\le \la^2.$$
It follows that
\begin{multline*}
\Bigl|\sum_{2\lambda<\lambda_j<\tau_k/2}\sum_{10\lambda<\tau_k< \la^2}\iint
\frac{\la_j^{2N}}{\la_j^2-\tau_k^2}\ala(\la_j)
 e^0_j(x)e^0_j(y)V(y)\tau_k^{-2N}
e_{\tau_k}(x)e_{\tau_k}(y) \, dxdy\Bigr|
\lesssim \la^{-\sigma},\ \forall \sigma.
\end{multline*}
So we just need to prove if $N=2n$, then
\begin{multline*}
\Bigl|\sum_{\lambda_j\le2\la}\sum_{10\lambda<\tau_k< \la^2}\iint
\frac{\la_j^{2N}}{\la_j^2-\tau_k^2}\ala(\la_j)
 e^0_j(x)e^0_j(y)V(y)\tau_k^{-2N}
e_{\tau_k}(x)e_{\tau_k}(y) \, dxdy\Bigr|
\lesssim \la^{n-2} .
\end{multline*}
 Let $$m_N(\la_j,s)=\frac{\la_j^{2N}}{1-s^2\la_j^2}\ala(\la_j).$$ Then for $s\in [0,(10\la)^{-1}]$ we have
\begin{equation}\label{mN}
|m_N(\la_j,s)|+|\la^{-1}\partial_s m_N(\la_j,s)|\ls \la^{2N}.\end{equation}
We can use the same argument as before to handle 
\[m_N(\la_j,\tau_k^{-1})=m_N(\la_j,0)+\int_0^{(2\la)^{-1}}\partial_s m_N(\la_j,s)\1_{[0,1/\tau_k]}(s)ds.\]
As before, we just need to handle the first term, and the second term is similar by \eqref{mN}. For $N=2n$ we have by \eqref{rough}
\begin{align*}
	&\Big|\sum_{\la_j\le 2\la}\sum_{10\la<\tau_k<\la^2}\iint\la_j^{2N}e_j^0(x)\ala(\la_j)e_j^0(y)V(y)\tau_k^{-2-2N}e_{\tau_k}(y)e_{\tau_k}(x)dxdy\Big|\\
	&\ls \|V\|_{L^1(M)}\cdot \la^{2N}\sup_y(\sum_{\la_j\le2\la}|e_j^0(y)|^2)^{1/2}(\sum_{\tau_k<\la^2}\tau_k^{-4-4N}|e_{\tau_k}(y)|^2)^{1/2}\\
	&\ls \|V\|_{L^1(M)}\cdot \la^{2N}\cdot\la^{n/2}\cdot \la^{\frac n2-2-2N}\\
	&\ls \la^{n-2}.
\end{align*}

$ $

\noindent \textbf{Case 2:} If $\lambda_j>2\tau_k$, then by the rapid decay property of $\ala(\lambda_j)$ we have
\begin{align*}
\Bigl| \, \sum_{\la_j>2\tau_k}
\sum_{\tau_k >10\lambda}\iint
\frac{\ala(\la_j)}{\la_j^2-\tau_k^2}
e_j^0(x)e_j^0(y)V(y)e_{\tau_k}(x)
e_{\tau_k}(y) \, dx dy\, \Bigr|
\lesssim \la^{-\sigma},\ \  \forall \, \sigma.
\end{align*}
As in Case 1, we only need to prove
\begin{equation}\label{high2}
\Bigl| \, \sum_{\la_j>2\tau_k}
\sum_{\tau_k >10\lambda}\iint
\frac{\ala(\tau_k)}{\la_j^2-\tau_k^2}
e_j^0(x)e_j^0(y)V(y)e_{\tau_k}(x)
e_{\tau_k}(y) \, dx dy\, \Bigr|\lesssim \la^{n-2} (\log\lambda)^\frac12.
\end{equation}
We similarly expand
$$
\frac{1}{\la_j^2-\tau_k^2}=\la_j^{-2}
+\la_j^{-2}\bigl(\tau_k/\la_j\bigr)^2+
\cdots + \la_j^{-2}\bigl(\tau_k/\la_j\bigr)^{2N-2}
\\
+(\tau_k/\la_j)^{2N}\frac{1}{\la_j^2-\tau_k^2}
$$
where we will choose $N=2n$ later. Then  we can repeat the argument in Case 1 (with $\la_j$ and $\tau_k$ interchanged) to obtain for $\ell=0,\dots,N-1$
\begin{align}\label{3.231}
\Bigl| \iint \sum_{\lambda_j>20\la}\la_j^{-2-2\ell} e_j^0(x)e_j^0(y)\sum_{10\la<\tau_k<\lambda_j/2}
\tau_k^{2\ell}
\ala(\tau_k) e_{\tau_k}(x)e_{\tau_k}(y)
\,
V(y) \, dxdy\Bigr|
\ls \la^{-\sigma}, \ \forall \sigma,
\end{align}
\begin{align}\label{3.241}
\Bigl|\sum_{\la_j>20\la}\sum_{10\la<\tau_k<\la_j/2}\iint
\frac{\tau_k^{2N}}{\la_j^2-\tau_k^2}
\ala(\tau_k) e^0_j(x)e^0_j(y)\la_j^{-2N}
e_{\tau_k}(x)e_{\tau_k}(y)V(y) dxdy\Bigr|
\ls \la^{-\sigma},\ \ \forall \sigma.
\end{align}
The bounds are better than \eqref{high2}, thanks to the rapid decay property of $\ala(\tau_k)$. So we complete the proof.

\bibliographystyle{plain}

\begin{thebibliography}{0}
\bibitem{ava}V. G. Avakumovi\'c. \"Uber die Eigenfunktionen auf geschlossenen Riemannschen Mannigfaltigkeiten.
Math. Z., 65:327–344, 1956.
\bibitem{berard}P. H. B\'erard. On the wave equation on a compact Riemannian manifold without conjugate points.
Math. Z., 155(3):249–276, 1977.
\bibitem{BSS}M. D. Blair, Y. Sire, and C. D. Sogge. Quasimode, eigenfunction and spectral projection bounds for
Schr\"odinger operators on manifolds with critically singular potentials. J Geom Anal 31, 6624–6661 (2021).

\bibitem{car1}T. Carleman, Propri\'et\'es asymptotiques des fonctions fondamentales des membranes vibrantes,
Comptes Rendus des Math\'ematiciens Scandinaves \'a Stockholm, (1934), pp. 14-18.
\bibitem{car2}T. Carleman, \"Uber die asymptotische Verteilung der Eigenwerte partielle Differentialgleichungen, Berichten der mathematisch-physisch Klasse der S\"achsischen Akad. der Wissenschaften zu Leipzig,
LXXXVIII Band, Sitsung, 15 (1936).
\bibitem{courant}R. Courant. \"Uber die Eigenwerte bei den Differentialgleichungen
der mathematischen Physik. Mat. Z., 7:1-57 (1920).
\bibitem{dan}Daners, Daniel. Heat kernel estimates for operators with boundary conditions. Math. Nachr. 217 (2000), 13–41.
\bibitem{davies}E.B. Davies. Gaussian upper bounds for the heat kernel of some second-order operators on Riemannian manifolds,J. Funct. Anal. 80(1988).16-32.
\bibitem{DG}Duistermaat, J.J., Guillemin, V.: The spectrum of positive elliptic operators and periodic bicharacteristics. Invent. Math. 29(1), 37–79 (1975)
\bibitem{fs}Frank, Rupert L.; Sabin, Julien.
Sharp Weyl laws with singular potentials. Pure Appl. Anal.5(2023), no.1, 85–144.
\bibitem{Gre}P. Greiner, An asymptotic expansion for the heat equation, Arch. Rat. Mech. Anal. 41 (1971), 168-218.
\bibitem{Gri}D. Grieser. Uniform bounds for eigenfunctions of the Laplacian on manifolds with boundary. Comm.
Partial Differential Equations, 27 (2002), 1283–1299.
\bibitem{Gu17}B. G\"uneysu. Heat kernels in the context of kato potentials on arbitrary manifolds. Potential Analysis, 46(1):119–134, 2017.
\bibitem{horm}L. H\"ormander. The spectral function of an elliptic operator. Acta Math., 121:193–218, 1968.
	\bibitem{hor4}L. H\"ormander. The Analysis of Linear Partial Differential Operators IV: Fourier Integral Operators, Springer-Verlag, Berlin, 1994.
\bibitem{hs}X. Huang and C.D. Sogge. Weyl formulae for Schr\" odinger operators with critically singular potentials.  Comm. Partial Differential Equations 46 (2021), no. 11, 2088–2133.
\bibitem{hz2021}X. Huang and C. Zhang, Pointwise Weyl Laws for Schrodinger operators with singular potentials. Adv. Math. 410 (2022), Paper No. 108688, 34 pp.

\bibitem{hz2023}X. Huang and C. Zhang, Sharp Pointwise Weyl Laws for Schrödinger Operators with Singular Potentials on Flat Tori, Comm. Math. Phys.401(2023), no.2, 1063–1125.
\bibitem{ivrii}Ivrii, V.: Second term of the spectral asymptotic expansion for the Laplace–Beltrami operator on manifold with boundary. Funct. Anal. Appl. 14(2), 98–106 (1980)
\bibitem{ivrii2}Ivrii, Victor. 100 years of Weyl's law. Bull. Math. Sci.6(2016), no.3, 379–452.
\bibitem{kac}M. Kac, Can one hear the shape of a drum? Amer. Math. Monthly 73, 1-23 (1966).
\bibitem{kac2}M. Kac, On some connections between probability theory and differential and integral equations,
Proc. Second Berkeley Symposium on Mathematical Statistics and Probability, 1950,
pp. 189-215, University of California Press, Berkeley and Los Angeles, 1951.
\bibitem{kst1}H. Koch, H. Smith and D. Tataru, Subcritical  $L^p$  bounds on spectral clusters for Lipschitz metrics. 
Math. Res. Lett. 15 (2008), no. 5, 993–1002.
\bibitem{kst2}H. Koch, H. Smith and D. Tataru, Sharp  $L^p$  bounds on spectral clusters for Lipschitz metrics. 
Amer. J. Math. 136 (2014), no. 6, 1629–1663.

\bibitem{lev}B. M. Levitan. On the asymptotic behavior of the spectral function of a self-adjoint differential
equation of the second order. Izvestiya Akad. Nauk SSSR. Ser. Mat., 16:325–352, 1952.
\bibitem{mel}R. Melrose, Weyl's conjecture for manifolds with concave boundary, Proc. SymposPure Math., vol.36,Amer. Math. Soc., Providence, RI, 1980, pp. 257-274.
\bibitem{MTbook}R. Melrose and M. Taylor. Boundary problems for the wave equation with grazing and gliding rays, manuscript
\bibitem{MS}McKean, H.P., JR., \& I. M. Singer, Curvature and the eigenvalues of the Laplacian. J.
Diff. Geometry 1, 43-69 (1967).

\bibitem{LY}Peter Li and Shing-Tung Yau, On the parabolic kernel of the Schr\"odinger operator, Acta Math. 156 (1986), no. 3-4, 153–201.
\bibitem{LHB}J. L\"orinczi, F. Hiroshima, and V. Betz. Feynman-Kac-type theorems and Gibbs measures on path
space: with applications to rigorous quantum field theory, volume 34. Walter de Gruyter, 2011.
\bibitem{RSS}G. V, Rozenblum, M. A. Shubin, and M. Z. Solomyak, Spectral theory of differentialoperators, Current problems in mathematics: fundamental directions, Partial Differential Equations VII, vol. 64, VINITI, Moscow, 1989; English transl., Encyclopaediaof Mathematical Sciences,vol, 64, Springer-Verlag, New York,1994.
\bibitem{safa}Safarov, Yu., Vassiliev, D.: The Asymptotic Distribution of Eigenvalues of Partial Differential Operators.
Translations of Mathematical Monographs. AMS, vol. 155 (1997)






\bibitem{seeley1}R. Seeley. A sharp asymptotic estimate for the eigenvalues of the
Laplacian in a domain of R3. Advances in Math., 102(3):244-264
(1978).
\bibitem{seeley2}R. Seeley. An estimate near the boundary for the spectral function
of the Laplace operator. Amer. J. Math., 102(3):869-902 (1980).
\bibitem{Simonbook}B. Simon. Functional integration and quantum physics, Academic Press, New York, 1979
\bibitem{simonsurvey}B. Simon. Schr\"odinger semigroups. Bull. Amer. Math. Soc. (N.S.), 7(3):447–526, 1982.
\bibitem{ssjams}H. Smith and C.D. Sogge. 
On the critical semilinear wave equation outside convex obstacles. J. Amer. Math. Soc.8(1995), no.4, 879–916.
\bibitem{ssacta}H. Smith and C.D. Sogge. On the  Lp  norm of spectral clusters for compact manifolds with boundary. 
Acta Math. 198 (2007), no. 1, 107–153.
\bibitem{soggebd}C.D. Sogge. 
Eigenfunction and Bochner Riesz estimates on manifolds with boundary. Math. Res. Lett.9(2002), no.2-3, 205–216.
\bibitem{hangzhou}C.D. Sogge. Hangzhou lectures on eigenfunctions of the Laplacian, volume 188 of Annals of Mathematics Studies. Princeton University Press, Princeton, NJ, 2014.
\bibitem{fio}C.D. Sogge, Fourier Integrals in Classical Analysis, Cambridge Tracts in Mathematics, vol.210,
Cambridge University Press, Cambridge, 2017.
\bibitem{sturm}Karl-Theodor Sturm, Schr\"odinger semigroups on manifolds, J. Funct. Anal. 118 (1993), no. 2, 309–350.

\bibitem{taylorpde1}Taylor, M.: Partial Differential Equations. I. Basic Theory. Applied Mathematical Sciences, vol. 115. Springer, New York (2010)
\bibitem{weyl}H. Weyl. \"Uber die Asymptotische Verteilung der Eigenwerte. Nachr.
Konigl. Ges. Wiss. G\"ottingen, 110-117 (1911).
\bibitem{weylconj}H. Weyl, \"Uber die Randwertaufgabe der Strahlungstheorie und asymptotische Spektralgesetze,
J. Reine Angew. Math. 143 (1913), 177–202.
\bibitem{xu}Xu, Xiangjin. Eigenfunction estimates for Neumann Laplacian and applications to multiplier problems.
Proc. Amer. Math. Soc. 139 (2011), no. 10, 3583–3599.

   \end{thebibliography}

\end{document}